# Optimization of Wireless Power Transfer Systems Enhanced by Passive Elements and Metasurfaces

Hans-Dieter Lang, *Student Member, IEEE,* and Costas D. Sarris, *Senior Member, IEEE*

*Abstract*—This paper presents a rigorous optimization technique for wireless power transfer (WPT) systems enhanced by passive elements, ranging from simple reflectors and intermediate relays all the way to general electromagnetic guiding and focusing structures, such as metasurfaces and metamaterials. At its core is a convex semidefinite relaxation formulation of the otherwise nonconvex optimization problem, of which tightness and optimality can be confirmed by a simple test of its solutions. The resulting method is rigorous, versatile, and general — it does not rely on any assumptions. As shown in various examples, it is able to efficiently and reliably optimize such WPT systems in order to find their physical limitations on performance, optimal operating parameters and inspect their working principles, even for a large number of active transmitters and passive elements.

## I. INTRODUCTION

W IRELESS power transfer (WPT) systems have been investigated in various forms, contexts and for various applications [1]–[4]. A large variety of unique solutions for powering and charging devices in an untethered fashion has been proposed; ranging from applications such as charging mobile devices [5] to electric vehicles [6]–[8] and including charging, operating as well as communicating with biomedical implants [9].

However, among other constraints, the performance — most commonly measured by the power transfer efficiency (PTE), the efficiency at which the power can be transferred from the transmitter(s) to the receiver — considerably limits the range of such applications. Usually, PTEs high enough for practical applications are confined to distances between transmitter(s) and receiver(s) in the order of fractions of the wavelength at the operating frequency as well as to dimensions comparable to those of the transmitter and receiver coils [5].

Various attempts to mitigate this problem and increase the PTE or extend the range of high PTE have been proposed. A popular method to enhance PTE involves introducing passive elements in the proximity of the transmitter(s) and the receiver [10], acting as relay elements or as parasitic elements that can enlarge the effective aperture of the transmitter.

Multiple passive elements can be used cooperatively; for example to form near-field guiding structures [11]–[15] or provide more general functions associated with metasurfaces or metamaterials [16]–[18].

As will be shown, the optimization of such systems is not a trivial task; the standard formulations of the required passivity and power constraints are nonconvex. This renders the whole problem unsolvable in general, particularly for a large number of unknowns (e.g. the number of passive elements). Here, however, a convex relaxation formulation is derived, which essentially constitutes a general and rigorous work-around to this problem. A simple test of the solution confirms tightness (meaning exact representation of the original problem) and hence, that the true global optimum has been found.

The resulting powerful and rigorous optimization method generalizes the previously presented optimization technique for WPT systems with multiple active transmitters [19]. Hence, it can be used to investigate the maximum performance of a particular system with multiple active transmitters and passive elements and to obtain the corresponding optimal operating parameters (i.e. currents and loading reactances). Furthermore, using an outer loop optimization, the optimum load resistance can also be found, leading to the maximum achievable PTE of the system; the absolute limit on the performance.

The outline of this paper is as follows: After a short introduction, the optimization problem of passively-enhanced WPT systems is presented. The convex semidefinite relaxation is then derived and a simple test for tightness is given. In the results section, three numerical examples prove the validity and versatility of the method: First, the case where additional passively excited transmitters are added behind the active transmitter is considered. Second, optimal relay configurations are treated, useful for both high-efficiency range extension as well as for misalignment mitigation. Lastly, the enhancement by general passive metasurfaces is discussed; proving the method's capability to handle a large number of passive elements. At the end, final remarks and conclusions are given.

Remarks on the notation: Thin italic letters represent scalar variables, bold small letters refer to vectors, bold capital letters are matrices; $\mathbf{v}^T$ stands for the transpose of the vector $\mathbf{v}$, while $\mathbf{v}^H$ stands for its Hermitian (conjugate transpose). The symbol $\succ$ ($\succeq$) is used to denominate positive (semi-) definiteness of matrices, respectively. A star $(\cdot)^\star$ marks the optimized arguments leading to the optimal solution, while a bar $(\bar{a})$ is used to denote the complex conjugate. Lastly, 'program' is used as a synonym for 'optimization problem', as is common in the context of mathematical optimization [20].

## II. PRELIMINARIES

### A. The MISO WPT System Model with Passive Elements

Fig. 1 depicts the general form of wireless power transfer (WPT) systems under consideration, incorporating multiple active transmitters, multiple passive elements and a single receiver. Systems with multiple transmitters will be referred to as MISO (multiple-input, single-output) WPT systems, in contrast to SISO (single-input single-output) systems with only a single transmitter.

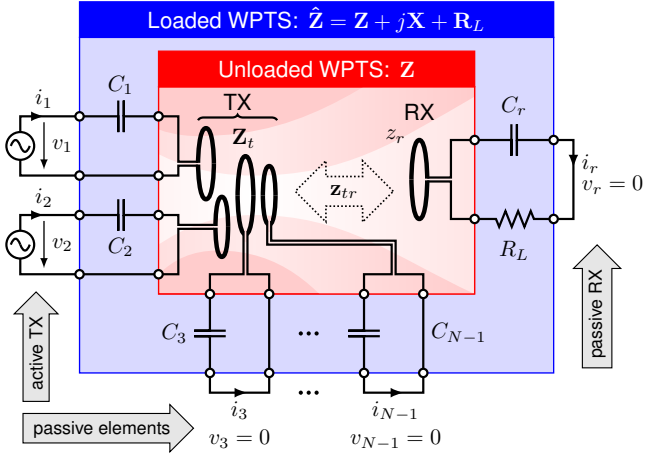

Fig. 1. Loop-based MISO WPT system model with multiple active transmitters and passive elements: Core structure with unloaded impedance matrix $\mathbf{Z}$, reactive components $x_n = -(\omega C_n)^{-1}$ (in case of capacitors), resistive receiver load $R_L$, and voltage sources $v_n$. In this example, subscripts 1 and 2 correspond to active transmitters (i.e. connected to a source), while 3 to $N-1$ stand for the passive elements.

Let the (unloaded) impedance matrix $\mathbf{Z} \in \mathbb{C}^{N \times N}$ of such systems be partitioned according to

$$\mathbf{Z} = \begin{bmatrix} \mathbf{Z}_t & \mathbf{z}_{tr} \\ \mathbf{z}_{tr}^T & z_r \end{bmatrix} = \begin{bmatrix} \mathbf{Z}_a & \mathbf{Z}_{ap} & \mathbf{z}_{ar} \\ \mathbf{Z}_{ap}^T & \mathbf{Z}_p & \mathbf{z}_{pr} \\ \mathbf{z}_{ar}^T & \mathbf{z}_{pr}^T & z_r \end{bmatrix} . \quad (1)$$

The subscripts $t$ and $r$ refer to the transmitter and receiver parts, respectively, whereas the subscripts $a$ and $p$ refer to the active and passive transmitters. Finally, the combinations $tr$, $ap$, $ar$ and $pr$ stand for the submatrices and vectors coupling two of these groups together.

Let all the nodes be numbered by $n = 1, ..., N$, where $N = A + P$ is the total number of nodes, including $A$ active transmitters, $n \in \mathcal{A}$ (note that $\mathbf{Z}_a \in \mathbb{C}^{A \times A}$) and $P$ passive nodes, $n \in \mathcal{P}$—note that this includes both the passive elements and the receiver.

The diagonals of $\mathbf{Z}$, $\mathbf{Z}_t$, $\mathbf{Z}_a$, and $\mathbf{Z}_p$ as well as $z_r$ refer to the loss resistances and self-reactances of the transmitters and receiver, respectively. Typically, all $Z''_{,nn}$, $z''_r > 0$ (inductive) when considering systems of magnetically coupled loops or coils; however, this is not a requirement. Likewise, all off-diagonal entries of $\mathbf{Z}$ are the coupling impedances due to the mutual inductances $j\omega M_{n,m}$. Generally each $M_{n,m}$ is complex (leading to non-zero real parts in the off-diagonal entries of $\mathbf{Z}$), due to retardation effects, when the electrical distance between the transmitter(s) and/or receiver are not very small.

For physical reasons, impedance matrices of passive reciprocal circuits have to be *positive-real* [21]: The impedance matrix is symmetric $\mathbf{Z} = \mathbf{Z}^T$ and all $\mathbf{Z}'$, $\mathbf{Z}'_t$, $\mathbf{Z}'_a$, $\mathbf{Z}'_p \succ 0$ and $z'_r > 0$ (positive-definite), meaning positive quadratic forms, for example:

$$\mathbf{i}^H \mathbf{Z}' \mathbf{i} > 0 \qquad \forall \, \mathbf{i} \in \mathbb{C}^N \quad (2)$$

Positive semidefiniteness ($\succeq 0$) implies nonnegative ($\geq 0$) quadratic forms. For negative (semi-) definiteness the opposite signs and directions apply. In the context of impedance matrices, the mathematical property of positive-definiteness represents the fact that there is always some amount of loss in the system.

The process of adding reactances to each of the transmitter and receiver nodes as well as adding a resistive load as receiver will be referred to as *loading* of the WPT system, where

$$\hat{\mathbf{Z}} = \mathbf{Z} + j\mathbf{X} + \mathbf{R}_L \quad (3)$$

is the resulting *loaded impedance matrix*, marked by a hat. In full detail, the voltages and currents $\mathbf{v}$, $\mathbf{i}$ of the entire WPT system are related according to

$$\underbrace{\begin{bmatrix} \mathbf{v}_a \\ \mathbf{v}_p = \mathbf{0} \\ v_r = 0 \end{bmatrix}}_{\mathbf{v}} = \left( \mathbf{Z} + j \underbrace{\begin{bmatrix} \mathbf{X}_a \\ & \mathbf{X}_p \\ & & x_r \end{bmatrix}}_{\mathbf{X}} + \underbrace{\begin{bmatrix} \mathbf{0} \\ & \mathbf{0} \\ & & R_L \end{bmatrix}}_{\mathbf{R}_L} \right) \underbrace{\begin{bmatrix} \mathbf{i}_a \\ \mathbf{i}_p \\ i_r \end{bmatrix}}_{\mathbf{i}} \quad (4)$$

where the subscripts $t$, $a$, $p$, and $r$ have the same meanings as explained for (1). The real-valued diagonal matrix $\mathbf{X}$ contains the load reactances, where positive and negative values refer to inductive and capacitive loading, respectively. They are grouped into the two diagonal matrices corresponding to the active transmitters and passive elements as well as a scalar for the receiver. The matrix $\mathbf{R}_L$ is zero everywhere but at the last diagonal entry, corresponding to the receiver, where the load resistance $R_L > 0$ is located.

It is important to note that, since the load resistance $R_L$ is actually part of the impedance matrix, KVL states that the corresponding receiver voltage is zero, i.e. $v_r = 0$ as explicated in (4) and depicted in Fig. 1. The same is true for the voltages of the passive couplers, $\mathbf{v}_p = \mathbf{0}$.

### B. The Objective: Power Transfer Efficiency (PTE)

The central figure of merit when optimizing WPT systems is the *power transfer efficiency*[1] (PTE) [1]–[5], [7], [9]: the ratio of the power $P_L$ transferred to the load $R_L$, to the total transmit power provided by the source $P_t$

$$\eta = \frac{P_L}{P_t} = \frac{P_L}{P_t + P_L} = \frac{R_L}{R_l + R_L} . \quad (5)$$

---

[1]It should be noted that the so defined PTE refers only to the electromagnetic power transfer efficiency; i.e. the impacts of impedance matching, rectification and generation onto the total efficiency of an entire wireless power transfer system are not considered, here.

$P_l$ is the total power absorbed by the system, due to dissipation and radiation, modeled by the loss resistance $R_l$.

The goal of this paper is to find a mathematical method to determine the absolute performance limits of such WPT systems to wirelessly transfer power most efficiently via multiple active transmitters and passive elements to a single receiver. Hence, the aim is to maximize the PTE, obtained according to its definition (5) as the biquadratic form

$$\eta = \frac{\frac{1}{2}\mathbf{i}^H \mathbf{R}_L \mathbf{i}}{\frac{1}{2}(\mathbf{i}^H \mathbf{v})'} = \frac{\mathbf{i}^H \mathbf{R}_L \mathbf{i}}{\mathbf{i}^H (\mathbf{Z}' + \mathbf{R}_L) \mathbf{i}} \ . \quad (6)$$

by finding the corresponding optimal voltages $\mathbf{v}$ and currents $\mathbf{i}$ (real and imaginary parts) as well as loading elements $\mathbf{x} = \mathrm{diag}(\mathbf{X}) = [\mathbf{x}_a, \mathbf{x}_p, x_r]$ and $R_L$. It is most convenient to formulate the PTE in terms of the currents $\mathbf{i}$. The voltages $\mathbf{v}$ can be obtained from the currents using the fully loaded system (4). However, both are only meaningful solutions as long as $v_r = 0$ and $\mathbf{v}_p = \mathbf{0}$.

## III. OPTIMIZATION

### A. Port Impedance Matrices (PIMs)

The total power inserted into the system $P_t$ in the denominator of (5) can be separated into the contributions of each transmitter as follows [19]:

$$P_t = P_l + P_L = \frac{1}{2}(\mathbf{i}^H \mathbf{v})' = \frac{1}{2}\mathbf{i}^H \hat{\mathbf{Z}}' \mathbf{i}$$
$$= \sum_n P_{t,n} = \frac{1}{2}\sum_n \mathbf{i}^H \hat{\mathbf{T}}_n \mathbf{i} \quad (7)$$

where $\hat{\mathbf{T}}_n$ are the (loaded) *port impedance matrices* (PIMs), and $n$ denotes the port. Note that $\mathbf{T}_n = \hat{\mathbf{T}}_n$, for all $n$ except $n = N$ (receiver), where $\hat{\mathbf{T}}_N = \mathbf{T}_N + \mathbf{R}_L$; i.e. loading only affects the $N$th PIM.

While the PIMs sum up to the total loss resistance matrices, i.e. $\sum_n \hat{\mathbf{T}}_n = \hat{\mathbf{Z}}' \succ 0$ (real-valued, symmetric and positive definite), the PIMs themselves are complex-valued, Hermitian, and indefinite. Furthermore, each PIM is singular, as it has exactly one positive, one negative and $N-2$ zero eigenvalues. These eigenvalues and eigenvectors can be obtained analytically, as shown in Appendix A.

### B. Quadratically-Constrained Quadratic Program (QCQP)

The goal is to minimize the total power required for unit transferred power $\frac{1}{2}\mathbf{i}^H \mathbf{R}_L \mathbf{i} = 1$ in terms of complex currents $\mathbf{i}$ under the following additional constraints:

- All power fluxes at ports of active nodes (transmitters) be nonnegative, i.e. $\frac{1}{2}\mathbf{i}^H \hat{\mathbf{T}}_n \mathbf{i} \geq 0$, $\forall n \in \mathcal{A}$ (the set of all $A$ active nodes), ensuring power is fed into the system via each active transmitter port, rather than drained from it [19].
- All power fluxes at ports of passive nodes (passive elements and the receiver) be zero, i.e. $\frac{1}{2}\mathbf{i}^H \hat{\mathbf{T}}_n \mathbf{i} = 0$, $\forall n \in \mathcal{P}$ (the set of all $P$ passive nodes). This represents the fact that there is neither power fed into nor drained from the system. Note that this still allows for power loss in that passive element.

Gathering all these constraints, the following program is obtained:

$$\begin{aligned}
P_{l,\min} = \min_{\mathbf{i}} \ & \frac{1}{2}\mathbf{i}^H \mathbf{Z}' \mathbf{i} && \text{(power loss)} \\
\text{s.t.} \ & \frac{1}{2}\mathbf{i}^H \hat{\mathbf{T}}_n \mathbf{i} \geq 0 && n \in \mathcal{A} \quad \text{(active nodes)} \\
& \frac{1}{2}\mathbf{i}^H \hat{\mathbf{T}}_n \mathbf{i} = 0 && n \in \mathcal{P} \quad \text{(passive nodes)} \\
& \frac{1}{2}\mathbf{i}^H \mathbf{R}_L \mathbf{i} = 1 && \text{(transferred power)}
\end{aligned} \quad (8)$$

This a nonconvex *quadratically constrained quadratic program* (QCQP) [20], because some (in this case all) constraints are quadratic and nonconvex. Such problems are well-known to be notoriously difficult to solve (as they belong to the class of NP-hard problems). Most importantly, there are no available algorithms which are able to solve these types of problems in general, for an arbitrary number of unknowns, in polynomial time.

### C. Tight Semidefinite Relaxation

The following derivation follows in essence the most common standard semidefinite relaxation techniques layed out in detail in any reference on the matter, see e.g. [22], and will therefore only be given in brief form.

*1) Nonconvex QCQP:* The phases of the currents are only defined relative to each other. Thus, the phase of one the currents can be chosen freely, without any loss in generality; any other absolute phase solution could then be obtained by renormalization. The most natural choice is a purely real receiver current, i.e. $i_r = i'_r$ and $i''_r = 0$. Thus, the actual problem size can be reduced to $M = 2N - 1$. The remaining unknowns are the currents (separated into their real and imaginary parts)

$$\mathbf{c} = \begin{bmatrix} \mathbf{i}'_t \\ i'_r \\ \mathbf{i}'' \end{bmatrix} = \begin{bmatrix} \mathbf{i}'_a \\ \mathbf{i}'_p \\ i'_r \\ \mathbf{i}''_a \\ \mathbf{i}''_p \end{bmatrix} \in \mathbb{R}^M \ . \quad (9)$$

Using the cyclic property of the trace, $\mathbf{c}^T \mathbf{Q}_i \mathbf{c} = \mathrm{tr}(\mathbf{c}^T \mathbf{Q}_i \mathbf{c}) = \mathrm{tr}(\mathbf{Q}_i \mathbf{c}\mathbf{c}^T)$, the objective and all constraints can be written in linear terms of the (quadratic) current matrix $\mathbf{C} = \mathbf{c}\mathbf{c}^T \succeq 0$:

$$\begin{aligned}
P_{l,\min} = \min_{\mathbf{c}} \ & \mathrm{tr}(\mathbf{Q}_0 \mathbf{C}) \\
\text{s.t.} \ & \mathrm{tr}(\mathbf{Q}_n \mathbf{C}) \geq 0 && n \in \mathcal{A} \\
& \mathrm{tr}(\mathbf{Q}_n \mathbf{C}) = 0 && n \in \mathcal{P} \\
& \mathrm{tr}(\mathbf{R}\mathbf{C}) = 1 \\
& \mathbf{C} = \mathbf{c}\mathbf{c}^T
\end{aligned} \quad (10)$$

where

$$\mathbf{Q}_n = \begin{bmatrix} \hat{\mathbf{T}}'_n & -\hat{\mathbf{T}}''_n \\ \hat{\mathbf{T}}''_n & \hat{\mathbf{T}}'_n \end{bmatrix}_{M \times M} \quad n = 1, \dots, N \tag{11}$$

$$\mathbf{Q}_0 = \begin{bmatrix} \mathbf{Z}' + \mathbf{R}_L & -\mathbf{Z}''_n \\ \mathbf{Z}''_n & \mathbf{Z}' + \mathbf{R}_L \end{bmatrix}_{M \times M} = \sum_{n=1}^{N} \mathbf{Q}_n \tag{12}$$

$$\mathbf{R} = \begin{bmatrix} \mathbf{R}_L & \mathbf{0} \\ \mathbf{0} & \mathbf{R}_L \end{bmatrix}_{M \times M} = \mathbf{Diag} \begin{bmatrix} \mathbf{0}^{N-1}, R_L, \mathbf{0}^{N-1} \end{bmatrix} \tag{13}$$

The subscripts of the matrices refer to the size of the leading principal submatrices thereof; the last row and column are dropped, because they would refer to the imaginary part of the receiver current, which is chosen to be zero at all times, as discussed.

Note that the programs (8) and (10) are mathematically fully equivalent. The nonconvexities of the quadratic inequality and equality constraints has been isolated in the last (nonconvex) equality constraint, requiring the matrix $\mathbf{C}$ to be the outer product of the actual vector of unknowns $\mathbf{c}$. All the other constraints are affine (or even linear) in $\mathbf{C}$.

*2) Semidefinite Relaxation:* As for any convex relaxation, the idea is to exchange the nonconvex constraint by a constraint that achieves almost the same — while at the same time being convex. The constraint $\mathbf{C} = \mathbf{c}\mathbf{c}^T$ is equivalent to the so-called *rank-1 condition* $\mathrm{rank}\,\mathbf{C} = 1$, which also implies $\mathbf{C} \succeq 0$. To obtain a convex semidefinite relaxation, the rank-1 condition is removed, while still requiring positive-semidefiniteness[2] .

Thus, the following *semidefinite program* (SDP) is obtained:

$$P_{l,\min}^{\text{relax}} = \min_{\mathbf{C}} \ \mathrm{tr}(\mathbf{Q}_0 \mathbf{C})$$
$$\text{s.t.} \ \mathrm{tr}(\mathbf{Q}_n \mathbf{C}) \geq 0 \quad n \in \mathcal{A}$$
$$\mathrm{tr}(\mathbf{Q}_n \mathbf{C}) = 0 \quad n \in \mathcal{P} \tag{14}$$
$$\mathrm{tr}(\mathbf{R}\mathbf{C}) = 1$$
$$\mathbf{C} \succeq 0 .$$

This is a convex program which can be solved reliably and efficiently using dedicated algorithms. In Matlab, it can be implemented comfortably using CVX [23], [24] and solved for example using the standard solver SDPT3 [25], [26].

*3) Duality:* The Lagrangian, the weighted (penalized) sum of the objective and constraints [20] for this SDP (14) is given by:

$$L = \mathrm{tr}\left( \left[ \mathbf{Q}_0 - \sum_{n \in \mathcal{A}} \lambda_n \mathbf{Q}_n - \sum_{n \in \mathcal{P}} \nu_n \mathbf{Q}_n - \sigma \mathbf{R} \right] \mathbf{C} \right) + \sigma \tag{15a}$$

$$= \mathrm{tr}\big( [\mathbf{P} - \sigma \mathbf{R}] \mathbf{C} \big) + \sigma \tag{15b}$$

$$= \mathrm{tr}\big( \mathbf{Q}\mathbf{C} \big) - \sigma \tag{15c}$$

where the dual variables (penalty factors) are:

[2] Another interpretation of $\mathbf{C} = \mathbf{c}\mathbf{c}^T$ is that $\mathbf{C} - \mathbf{c}\mathbf{c}^T = 0$, equivalent to the two semidefinite constraints $\mathbf{C} - \mathbf{c}\mathbf{c}^T \succeq 0$ *and* $\mathbf{C} - \mathbf{c}\mathbf{c}^T \preceq 0$. Then, the relaxation takes place by removing the second constraint and using the first alone, by requiring the matrix of which $\mathbf{C} - \mathbf{c}\mathbf{c}^T$ is the Schur complement, to be positive-semidefinite [19].

- $\lambda_n \geq 0$, in order to only punish $\mathrm{tr}(\mathbf{Q}_n \mathbf{C}) < 0$, $\forall n \in \mathcal{A}$,
- $\nu_n \in \mathbb{R}$, penalizing $\mathrm{tr}(\mathbf{Q}_n \mathbf{C}) \neq 0$, $\forall n \in \mathcal{P}$, and finally
- $\sigma \in \mathbb{R}$, responsible for punishing $\mathrm{tr}(\mathbf{R}\mathbf{C}) \neq 1$.

To shorten the notation, the first two groups of dual variables are referred to as vectors $\boldsymbol{\lambda} = [\lambda_n]$ and $\boldsymbol{\nu} = [\nu_n]$, respectively.

The *primal problem* (i.e. the original problem (14)) is obtained from the Lagrangian by minimizing the objective while an inner maximization ensures that all constraints are satisfied:

$$P_{l,\min}^{\text{primal}} = \min_{\mathbf{C} \succeq 0} \max_{\substack{\boldsymbol{\lambda} \succeq 0 \\ \boldsymbol{\nu}, \sigma}} L . \tag{16}$$

Note that, if one or more of the constraints were not satisfied, the Lagrangian would be unbounded and diverge to $+\infty$. Thus, the primal problem is (16) is really identical to (14).

The so-called *dual problem* is obtained by interchanging the order of the maximization and minimization:

$$P_{l,\min}^{\text{dual}} = \max_{\substack{\boldsymbol{\lambda} \succeq 0 \\ \boldsymbol{\nu}, \sigma}} \min_{\mathbf{C} \succeq 0} L = \max_{\substack{\boldsymbol{\lambda} \succeq 0 \\ \boldsymbol{\nu}, \sigma}} \{ \sigma \ : \ \mathbf{Q} \succeq 0 \} . \tag{17}$$

Note that the dual program (17) takes the same form as in [19]. However, some constraints on the dual variables $\lambda_n \geq 0$ (active transmitters) within $\mathbf{Q}$ have been relaxed to $\nu_n \in \mathbb{R}$ (passive elements), thereby broadening the solution space and allowing for larger $\sigma$, meaning lower optimal PTE $\eta$.

As long as the primal problem satisfies certain conditions [20] (which are satisfied in this case, since the problem is continuous and convex), the dual optimum (17) is equivalent to the primal optimum (16). Therefore, in this case, the optimal solutions of the two problems are the same, $P_{l,\min}^{\text{primal}} = P_{l,\min}^{\text{dual}}$.

*4) Optimality conditions:* The Karush-Kuhn-Tucker (KKT) conditions of optimality [20] for the two problems are:

$$\mathbf{C}^\star \succeq 0 \quad \text{(primal feasibility)} \tag{18a}$$
$$\mathbf{Q}^\star \succeq 0 \quad \text{(dual feasibility 1)} \tag{18b}$$
$$\boldsymbol{\lambda}^\star \succeq 0 \quad \text{(dual feasibility 2)} \tag{18c}$$
$$\mathrm{tr}(\mathbf{R}\mathbf{C}^\star) = 1 \quad \text{(primal equality 1)} \tag{18d}$$
$$\mathrm{tr}(\mathbf{Q}_n \mathbf{C}^\star) \geq 0 \quad n \in \mathcal{A} \quad \text{(primal inequalities)} \tag{18e}$$
$$\mathrm{tr}(\mathbf{Q}_n \mathbf{C}^\star) = 0 \quad n \in \mathcal{P} \quad \text{(primal equalities)} \tag{18f}$$
$$\lambda_n \mathrm{tr}(\mathbf{Q}_n \mathbf{C}^\star) = 0 \quad n \in \mathcal{A} \quad \text{(compl. slackness 1)} \tag{18g}$$
$$\mathrm{tr}(\mathbf{Q}^\star \mathbf{C}^\star) = 0 \quad \text{(complementary slackness 2)} \tag{18h}$$

The last condition, the complementary slackness, boils down to $\mathbf{Q}^\star \mathbf{C}^\star = \mathbf{0}$, since $\mathbf{C}^\star, \mathbf{Q}^\star \succeq 0$. In the case of $\mathrm{rank}\,\mathbf{C}^\star = 1$ this also means that $\mathrm{rank}\,\mathbf{Q}^\star = M - 1$, or in other words that $\mathbf{Q}^\star$ has exactly one zero eigenvalue and $\mathbf{c}^\star$ is the corresponding eigenvector, the sought-after solution to the problem

*5) Tightness:* Since the constraints of the actual problem (10) have been loosened (essentially simply by removing one of them), the optimal value of (14) is generally only a lower bound on the actual optimum, i.e. $P_{l,\min}^{\text{relax}} \leq P_{l,\min}$. Therefore, by application of (5), the solution is only an upper bound on the maximum possible PTE, i.e. $\eta_{\max}^{\text{relax}} \geq \eta_{\max}$.

Equality is only achieved in the case of (14) being a *tight* (meaning exact) relaxation of (8) and (10). This would require that the removed rank constraint is naturally satisfied by the rest of the problem.

To test tightness, either the eigenvalues of $\mathbf{C}^\star$ or $\mathbf{Q}^\star$ could be investigated, or, numerically more efficient, the (normalized) *tightness error*, defined as

$$\epsilon = \frac{\left\| \mathbf{C}^\star - \mathbf{c}^\star (\mathbf{c}^\star)^T \right\|}{(\mathbf{c}^\star)^T \mathbf{c}^\star} \tag{19}$$

can be used. Whenever $\epsilon$ is small (i.e. within the bounds of typical numerical approximation errors of such algorithms), the semidefinite relaxation is tight and provides an exact optimum solution $\mathbf{c}^\star$ to the nonconvex QCQP (10).

Numerical experiments showed that this seems always the case for practical setups (typically down to $\epsilon \approx 10^{-14...-10}$). Thus, $\mathbf{C}^\star$ has exactly one nonzero (and indeed positive) eigenvalue and $\mathbf{c}^\star$ is its corresponding eigenvector; the global optimum to (10).

### D. Optimal Load Reactances and Voltages

The reactive loads in all passive nodes (including the receiver) are obtained from solving the corresponding homogeneous linear equations $\mathbf{v}_p = 0$ and $v_r = 0$. Physically, this means that in each of these loops, the voltages across these reactive elements are the same but with opposite sign as across everything else (the unloaded WPT system), in order to net in zero voltages overall.

To this end, the unloaded voltages are calculated according to

$$\mathbf{v}^\star_{\text{unloaded}} = \mathbf{Z}\mathbf{i}^\star . \tag{20}$$

Then, the optimal load reactances for the passive nodes ($n \in \mathcal{P}$) can be obtained using

$$x^\star_n = -\left( \frac{v^\star_{\text{unloaded},n}}{i^\star_n} \right)'' . \tag{21}$$

It can easily be verified that the zero real power flux constraint at port $n$, $P_n = \frac{1}{2}(\mathbf{i}^\star)^H \hat{\mathbf{T}}_n \mathbf{i}^\star = \frac{1}{2}(v^\star_{\text{unloaded},n} \bar{i}^\star_n)' = 0$, leads to $x^\star_n i^\star_n = -v^\star_{\text{unloaded},n}$ as required, and, thus, to zero voltage $v_n = 0$ at that port in the loaded case, when applying (4).

### E. Optimization of the Load Resistance

The presented optimization method provides optimal currents and reactive loading components, for a particular load resistance $R_L$. This is a very practical case; as the load resistance is usually not free to be chosen. However, in some cases, particularly when investigating the maximum achievable performance of the WPT in general, also the optimal load resistance $R_L$ is of importance.

To this end, the optimal load resistance can be obtained in an outer optimization loop, as follows:

$$R^\star_L = \arg\min_{R_L} P^{\text{relax}}_{l,\min}(R_L) \tag{22}$$

Numerical experiments confirm the expectation, based on practical experience with dissipation of physical multiport systems in general, that this remaining outer optimization of finding the optimal load resistance is always convex.

Note that the optimal performance is generally not very sensitive to the load resistance [27]. Thus, for practical considerations, the closed-form $R^\star_L$ for optimal fully active MISO configurations provided from the framework presented in [4] is usually very close to the true optimum and a simple and computationally efficient alternative. The calculation only requires the *minimum-loss output impedance* $z_o$ and (the square of) the *mutual coupling quality factor* $U$, both directly obtainable from the unloaded impedance matrix:

$$z_o = z_r - \mathbf{z}_{tr} \mathbf{Z}_t^{-1} \mathbf{z}'_{tr} \tag{23}$$

$$U^2 = \frac{\mathbf{z}_{tr}^H \mathbf{Z}_t^{-1} \mathbf{z}_{tr}}{z'_o} \tag{24}$$

$$R^\star_L = z'_o \sqrt{1 + U^2} \tag{25}$$

where the partitions are according to (1). As discussed in [27], these are the physically meaningful generalizations of the well-known SISO WPT systems [3].

## IV. APPLICATION EXAMPLES

It is well understood that the region of high PTE is usually limited to electrically very short distances between transmitter(s) and receiver [5]. Several approaches have been proposed in order to extend the electrical diameter of the high-PTE region, leading to the so-called "mid-range" [9] wireless power transfer systems. The most straightforward ideas are shielding, guiding (near-field beamforming) and focusing; all of which will be addressed briefly in the following via the three examples depicted in Fig. 2.

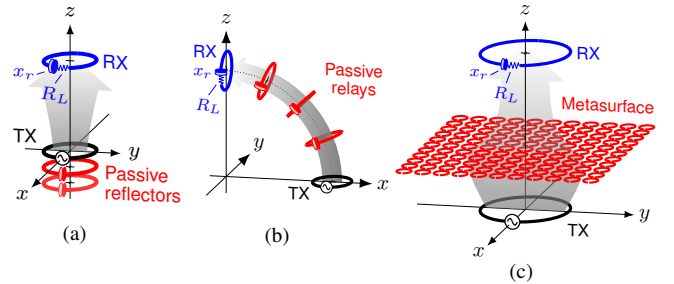

Fig. 2. The three principles of WPT performance enhancement under consideration: (a) Shielding using passive elements as reflectors, (b) guiding the magnetic near-field with passive relays [11], and (c) most general enhancement by metasurface (or -material), e.g., as proposed in [16], [17].

Note: The aim is not to thoroughly investigate each enhancement technique, but to show the versatility of the aforementioned optimization technique via those three relevant numerical examples.

### A. Enhancement via Shielding by Passive Loops

One of the most straightforward ideas for WPT performance enhancement is to try to confine the magnetic fields to the regions between the transmitter(s) and the receiver. This is often accomplished using a ferrite shield backing the transmitter(s) and receiver(s) in directions away from each other [8]. However, due to the lossy nature of ferrites, their application is usually restricted to the very low frequency regions. For mobile applications also their bulkiness, weight and cost (or availability) can be an issue. Recently, using resonant structures instead has been proposed, to reduce field leakage in unwanted directions [6].

The proposed test setups as shown in Fig. 3 are adopted from active MISO WPT system investigations [4]: one and two additional transmitter loops, respectively, are placed behind the first transmitter (thus at greater distance to the receiver), at a constant separation of $\Delta z = \lambda/1000$, where $\lambda$ is the freespace wavelength of the operating frequency $f = 13.56\,\text{MHz}$.

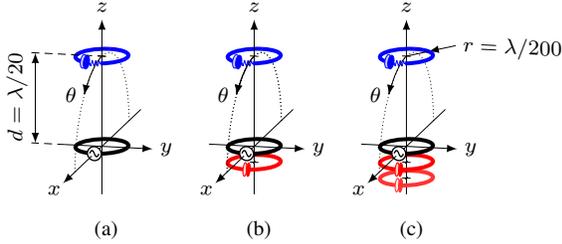

Fig. 3. The standard (reference) SISO WPT system (a) and two cases with appended passive elements: (b) the SISO-P1 with one passive and (c) the SISO-P2 with two passive elements, separated at $\Delta z = \lambda/1000$. Operating frequency $f = 13.56\,\text{MHz}$, wire thickness $t = 5\,\text{mm}$.

It has been demonstrated for MISO WPT systems that such constellations can outperform the standard SISO setup (a), when all transmitters are active and excited optimally. However, in this case, the additional loops are used as passive elements only, meaning they do not contain any active excitation, but are excited via passive coupling only. The passive loops are shorted by (lumped) reactive elements whose reactances are then optimized using the proposed method, along with the receiver reactance, to obtain the optimally enhanced power transfer from the transmitter (black) to the receiver (blue). As mentioned, the load resistance $R_L$ in the receiver could also be optimized, but since the impact on the final performance is usually negligible, it is omitted, here. The SISO WPT setups with one and two passive elements are called SISO-P1 and SISO-P2, respectively.

Fig. 4 compares the resulting maximum PTEs of the active MISO solutions (with multiple active transmitters) and the passive configurations in all three cases, along with the tightness errors. Interestingly, the maximum achievable PTEs are very close (a); even relative differences (b) are very small and only noticeable at angles where power transfer is difficult in general. As can be seen, similar to the active MISO cases, adding one or two passive loops increases the

PTE from about 13% to about 21% almost 26%, respectively. Thus, adding passive loops — even at a greater distance to the receivers — can lead to substantial enhancement of the performance. However, when comparing the currents and power contributions of each of the loops, as listed in Tab. I, it is observed that the optimal operating conditions for these additional loops lead to currents comparable to the MISO cases. Thus, when operated optimally, the additional loops do not act as reflectors or shielding. Instead, they seem to be used as additional transmitters, excited passively from the remaining active transmitter, but otherwise very similar to the active MISO case.

Lastly, as Fig. 4(c) confirms, apart from areas of general numerical difficulties in the angles around $\pm 60^\circ$, the tightness errors are in the range of $10^{-14}$ to $10^{-12}$. Thus, the relaxation is tight, meaning (14) is equivalent to the original (nonconvex) program (8) and the one true global optimum has been found.

TABLE I
CURRENTS AND POWER CONTRIBUTIONS OF THE ACTIVE MISO-2/-3 SETUPS COMPARED TO THE PASSIVELY ENHANCED SISO-P1/-P2 CASES.

| $n$ | **MISO-2** (two active TX) | | **SISO-P1** (one passive elem.) | |
|---|---|---|---|---|
| | $i_n$ (in A) | $P_n$ (in W) | $i_n$ (in A) | $P_n$ (in W) |
| 1 | $22.6 \angle 91.3^\circ$ | 4.60 | $22.6 \angle 91.4^\circ$ | 12.7 |
| 2 | $22.2 \angle 91.3^\circ$ | 8.12 | $22.2 \angle 91.3^\circ$ | 0 |
| Total | ($\eta = 7.86101\%$) | 12.7 | ($\eta = 7.86100\%$) | 12.7 |

(a)

| $n$ | **MISO-3** (three active TX) | | **SISO-P2** (two passive elem.) | |
|---|---|---|---|---|
| | $i_n$ (in A) | $P_n$ (in W) | $i_n$ (in A) | $P_n$ (in W) |
| 1 | $15.4 \angle 91.4^\circ$ | 1.07 | $15.4 \angle 91.5^\circ$ | 9.28 |
| 2 | $15.1 \angle 91.4^\circ$ | 3.03 | $15.1 \angle 91.4^\circ$ | 0 |
| 3 | $14.8 \angle 91.4^\circ$ | 5.18 | $14.8 \angle 91.4^\circ$ | 0 |
| Total | ($\eta = 10.77904\%$) | 9.28 | ($\eta = 10.77903\%$) | 9.28 |

(b)

### B. Enhancement via Passive Relays

Another intuitive method to enhance WPT over larger distances is placing a passive element in-between the transmitter and the receiver, so that it may act as a relay, over which the power transfer takes place [10]. With multiple such relays entire relay-line-type guiding structures can be constructed [11], [12] and could potentially even facilitate power distribution over multiple paths and to multiple loads [13]–[15]. In addition to extending the high-efficiency range of the power transfer, relays can also help mitigate issues due to misalignment between the transmitter and the receiver.

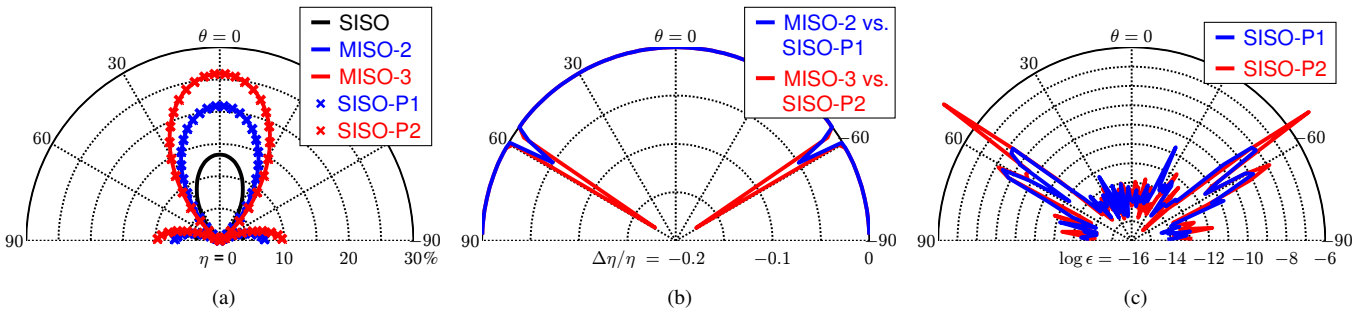

Fig. 4. Performance comparisons (a) of the full (active) MISO (lines) and the SISO cases with passive elements (x markers): MISO-2 vs. SISO-P1 (blue) and MISO-3 vs. SISO-P2 (red). The difference (b) is often almost unnoticeable. The tightness errors $\epsilon$ in (c) confirm that the relaxations in for both setups (blue SISO-P1, red SISO-P2) are equivalent to the original problems and that the global optima have been found.

Both possibilities are investigated in the following example, starting off from the SISO setup as shown in Fig. 5(a): The transmitter and receiver (radius $r_{tx} = r_{rx} = \lambda/200$) are placed at either ends of a quadrant of a circle of radius $r = \lambda/20$, at the operating frequency $f = 13.56\,\text{MHz}$. Even though they are placed at an electrically small distance ($d \approx \lambda/14$), the maximum achievable PTE of this system is only about 13% (b) — due to their severe misalignment.

As passive relay elements are inserted into the circular path as shown in Fig. 5(c)-(h), the situation changes considerably; higher optimum load resistances result and the overall performance of the WPT system is enhanced considerably. As can be seen, a single relay (c), leads to the largest performance enhancement step: the PTE increases from 13.2% to over 60%. Interestingly, this is only about 9% under what could have been achieved when using the passive relay as active transmitter. Thus, the passive relay setup is performing considerably better than a series of two SISO WPT systems of half the distance: $(69.5\%)^2 \approx 47.9\% < 60.2\%$. Moreover, using the relay, the optimum load resistance more than quadrupled. Increasing the number of relays leads to continued performance enhancement, but the increase in PTE slows down. Nevertheless, four relays (f) lead to almost 90% PTE, ten relays are able to exceed 95%. As the field plots reveal, with adding relays, the magnetic field is more and more confined to the inside of this guiding structure, forming a bent beam from the transmitter to the receiver and reducing the overall outside fields considerably, while the transferred power remains constant.

The graph in Fig. 6(a) shows the optimal performances vs. number of relays for different radii of the circular paths. The previously mentioned rapid improvement of the performance when using only a small number of relays is only observed for small radii; at greater distances a considerably larger number of relays is needed to achieve substantial performance enhancement. However, even there, eventually the PTE increases substantially and seems to asymptotically approach a value close to 100%, only restrained by conduction loss of the relays. As expected, unloaded (open-connected) relays remain invisible (solid black line), whereas shorted relays would lead to decreased PTE (dotted black line).

The graphs in Fig. 6(b) again prove that the relaxations were tight in all cases (with tightness errors in the range of $10^{-14}$ to $10^{-10}$) and that, therefore, the true global optima were attained.

Whereas the optimization of all reactive loading elements

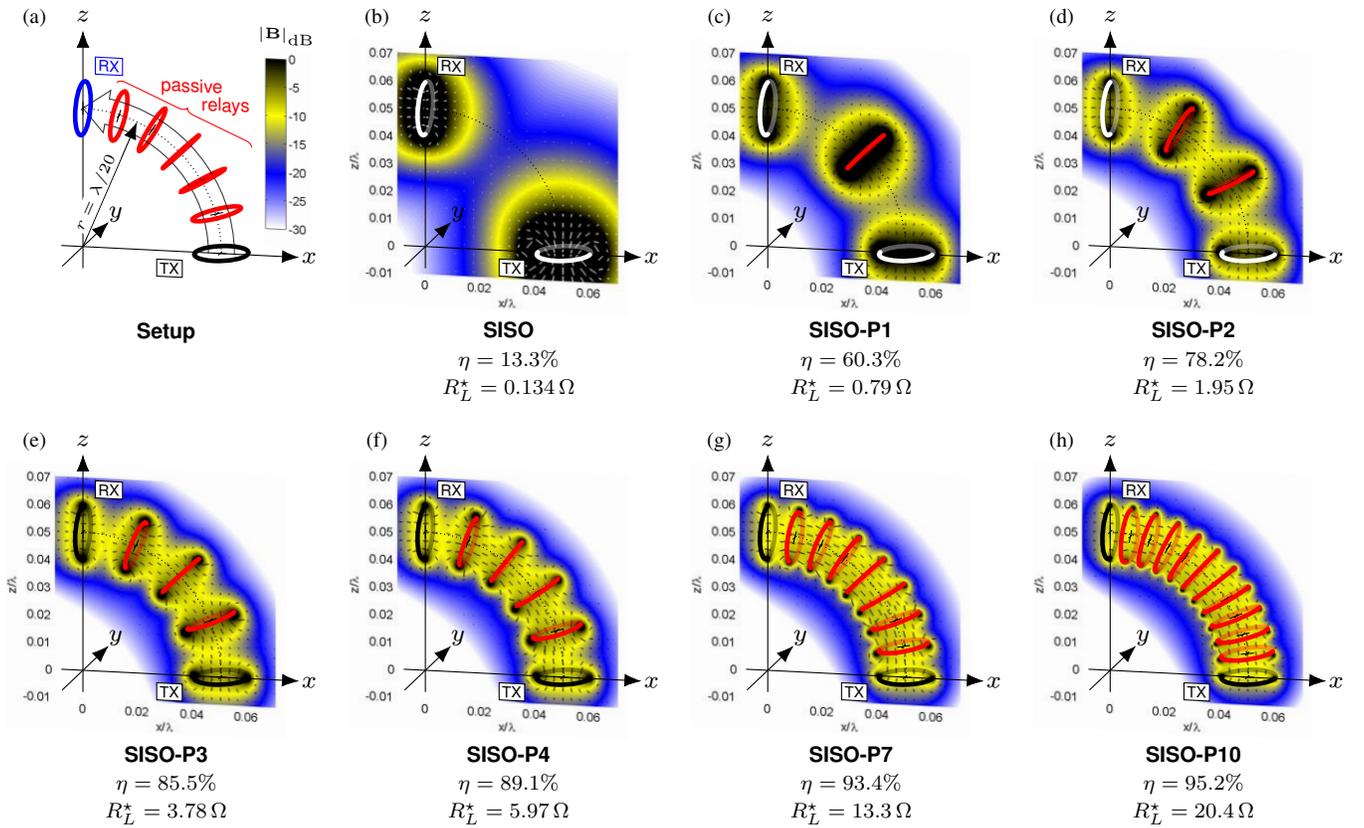

Fig. 5. Comparison of some relay-enhanced WPT configurations, guiding power transfer along an arc (of radius $r = \lambda/20$). Adding relays enhances the WPT performance (max. achievable PTE) not only by constraining the fields along the arc, but also by mitigating the issue of misalignment. The performance increase of adding a single relay is the largest; inserting more relays improves the overall performance further and also continues to lead to larger optimal load resistances.

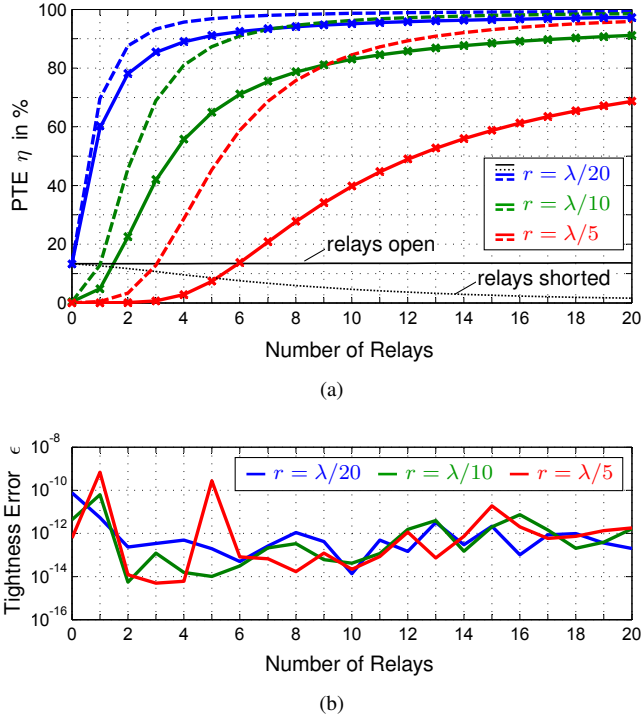

(a)

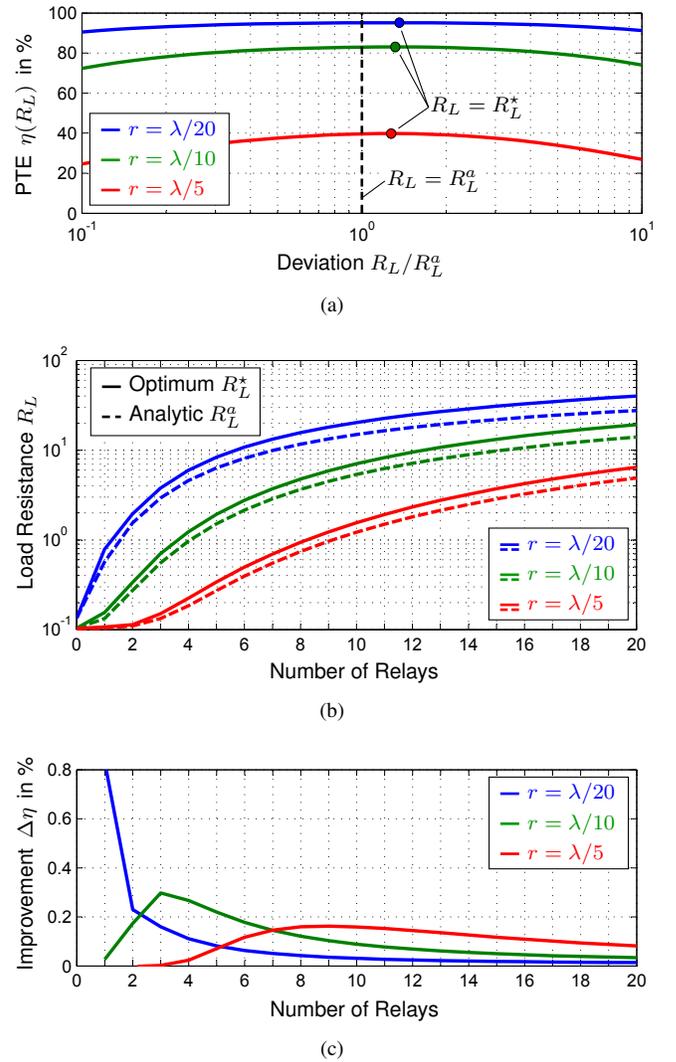

(b)

Fig. 6. Comparison of the performance enhancement due to optimal intermediate passive elements ("relays") for three distances (solid lines) and the maximum PTEs if the closest relay would be used as active transmitter, serving as an upper bound (dashed). This bound can be almost arbitrarily approached with increasing number of relays (a). The relaxations are tight in all cases (b).

(in these cases capacitive) of the passive relays is critical for achieving optimal performance, it is not sensitive to the receiver load resistance. This is addressed in Tab. II and 7, where the PTEs and loading capacitances are compared for cases with the optimum load resistance $R_L^\star$ and with estimates obtained from the closed-form expressions for active MISO WPT systems $R_L^a$. It can be seen that the difference in performance is small, even though the differences in the used load resistances can be rather large. The examples in Fig. 7(a) underline the fact that the optimum PTE is rather robust with respect to the actual $R_L$ around the optimum value $R_L^\star$ and that the PTE is a well-behaved concave function of the load resistance $R_L$, enabling numerical optimization

(a)

(b)

(c)

Fig. 7. Assessment of the importance of optimizing the receiver resistance $R_L$ for achieving optimal performance: The maximum PTE is a well-behaved concave function of $R_L$ (a), well-suited for optimization from the closed-form estimate $R_L^a$. However, even though the difference in load resistance can be significant (b), the difference in performance is usually negligible (c).

in an outer loop around the presented semidefinite relaxation with $R_L^a$ as starting point. While the difference between the optimum load $R_L^\star$ and the estimate $R_L^a$ grows larger with an increasing number of relays, as shown in Fig. 7(b), the resulting difference in PTE is decreasing, plotted in Fig. 7(c). The resulting capacitances are in a practically realistic region of around a hundred picofarads, with the relay capacitances depending more heavily on the number of relays than the receiver reactance. Furthermore, in the case with optimum receiver resistance $R_L^\star$, the resulting capacitance values appear to always be symmetric along the guiding path.

TABLE II
Comparison of performances and reactive loads for SISO setups with up to five passive relays when using the optimum receiver resistance vs. using estimates obtained from the active MISO closed-form expressions.

| Setup | | SISO | | SISO-P1 | | SISO-P2 | | SISO-P3 | | SISO-P4 | | SISO-P5 | |
|---|---|---|---|---|---|---|---|---|---|---|---|---|---|
| $R_L =$ | | $R_L^\star$ | $R_L^a$ | $R_L^\star$ | $R_L^a$ | $R_L^\star$ | $R_L^a$ | $R_L^\star$ | $R_L^a$ | $R_L^\star$ | $R_L^a$ | $R_L^\star$ | $R_L^a$ |
| $\eta_{\max}$ | (%) | 13.3 | 13.3 | 60.3 | 59.5 | 78.2 | 78.0 | 85.6 | 85.3 | 89.1 | 89.0 | 91.1 | 91.1 |
| $R_L$ | ($\Omega$) | 0.134 | 0.134 | 0.79 | 0.57 | 1.95 | 1.56 | 3.78 | 2.97 | 5.97 | 4.66 | 8.35 | 6.49 |
| $C_r^\star$ | (pF) | 104.43 | 104.43 | 104.38 | 104.32 | 104.35 | 104.29 | 104.20 | 104.02 | 103.98 | 103.69 | 103.68 | 103.27 |
| $C_p^\star$ | (pF) | | | 104.29 | 104.22 | 103.35 | 103.55 | 102.53 | 102.51 | 101.11 | 101.10 | 99.40 | 99.40 |
| | | | | | | 103.35 | 103.20 | 100.80 | 101.28 | 98.79 | 98.82 | 96.16 | 96.15 |
| | | | | | | | | 102.53 | 102.11 | 98.79 | 99.48 | 96.26 | 96.38 |
| | | | | | | | | | | 101.11 | 100.31 | 96.16 | 97.01 |
| | | | | | | | | | | | | 99.40 | 98.14 |

### C. Enhancement via General Metasurfaces and -Materials

Finally, the most general form of enhancement of wireless power transfer is considered: Enhancement by metamaterials and metasurfaces. Such enhancements have been considered by many research groups, see for example [16]–[18]. For simplicity and brevity, only an example with a metasurface (MS) will be considered; the generalization to (bulk-) metamaterials is straightforward. Instead of using any direct methods (e.g. Huygens surface approaches using the equivalence principle) to determine the surface impedance and obtain the required metasurface elements from that, the approach here is to optimize the loading of all of the elements of a particular metasurface with the semidefinite-relaxation based method presented herein.

Fig. 8 shows the setup of the metasurface-enhanced WPT example under consideration: The transmitter and receiver loops (of diameters $D_{\text{tx}} = D_{\text{rx}} = \lambda/100$) are separated by distance $d$. In the middle, a $15 \times 15$-element planar loop array (loop diameters $\lambda/500$, arranged periodically at a uniform pitch of $\lambda/200$ in the $x$- and $y$-directions) is inserted, which is referred to as metasurface. Each of the loops is loaded by a reactance (typically a capacitor $C_{n,m}$), which is to be optimized for maximum power transfer efficiency together with the receiver reactance $x_r$. The closed-form expression (25) is used to obtain the receiver resistance $R_L$.

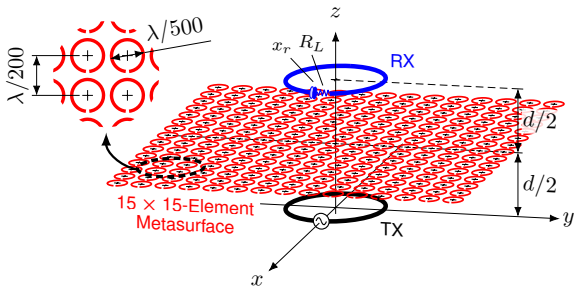

Fig. 8. Setup of the WPT system enhanced by a passive $15 \times 15$-element loop-based metasurface (MS). Dimensions: transmitter and receiver loop diameters $D_{\text{tx}} = D_{\text{rx}} = \lambda/100$ ($\approx 44\,\text{cm}$), wire thicknesses $\lambda/2000$ ($\approx 2.2\,\text{cm}$) for transmitter and receiver and $\lambda/10^4$ ($\approx 4.4\,\text{mm}$) for metasurface loops, respectively. Distance $d = 0.01 \ldots 0.1\lambda$ ($\approx 0.44 \ldots 4.4\,\text{m}$) at the operating frequency of $f = 6.78\,\text{MHz}$.

This case can be optimized in straightforward manner with the presented method. The 225 elements of the $15 \times 15$ loop surface in this case lead to a problem size of over 200000 unknowns, since the quadratic current matrix $\mathbf{C}$ has $(2 \times (225 + 2) - 1)^2$ entries. However, the problem is still solvable within a moderate computation time of a few minutes (on an Intel i3 machine running CVX/SDPT3 in Matlab).

The graphs in Fig. 9(a) present the fundamental results. As reference and lower bound, the black dotted line and the markers denote the cases without the metasurface (all elements open-connected) or when all metasurface loops are shorted; since the loops are electrically very small, the performances of both of these cases are very similar. The most important result is the solid red line, giving the maximum achievable performance with an optimum passive metasurface present. For comparison, the dashed red line plots the PTE of an active metasurface, where each loop of the metasurface is connected to a generator, serving as an upper bound on the achievable performance.

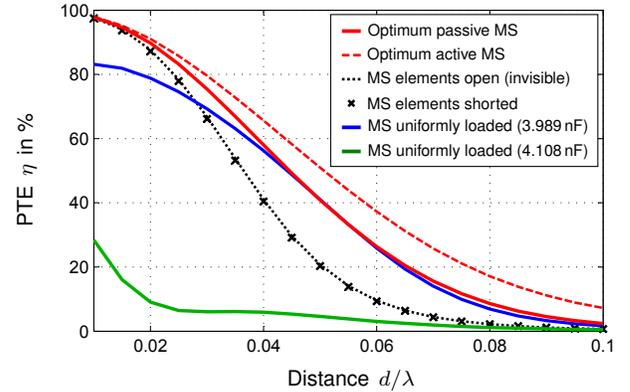

(a)

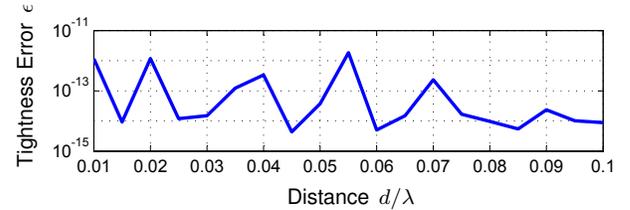

(b)

Fig. 9. Comparison of PTEs of the metasurface-enhanced WPT system: (a) Optimum active and passive PTEs (dashed and solid red, respectively), open and shorted metasurface loops (black), and two cases of uniform loading: a "good average" (blue) vs. the worst-case (green). The tightness error of the optimum passive metasurface (b) proves tightness of the relaxation.

As can be seen, the enhancement due to the passive surface is significant: For example at distance $d = 0.05\lambda$, the passive surface leads to an additional 20% in PTE over the reference case without a metasurface. On the other hand, at the same distance, the performance with the passive metasurface is only about 10% behind that with an active metasurface.

The blue and green solid curves in Fig. 9(a) address the importance of adequate loading of the metasurface elements. In each of these two cases, all loops of the metasurface are loaded uniformly by the same capacitance: about 3.989 nF (the optimal center capacitance at $d = 0.06\lambda$) in the blue case and 4.108 nF in the green case. The blue curve represents a case of good overall performance, which after some distance (about $0.04\lambda$) is largely comparable to the optimum (solid red). The green curve represents a worst-case scenario: Using an only slightly different capacitance value (differing by about 3%), degrades the maximum overall performance far below the case without a metasurface (the receiver reactance is optimized to still maximize the performance using the suboptimally loaded metasurface). Hence, the loading of the metasurface is absolutely crucial for optimal performance enhancement.

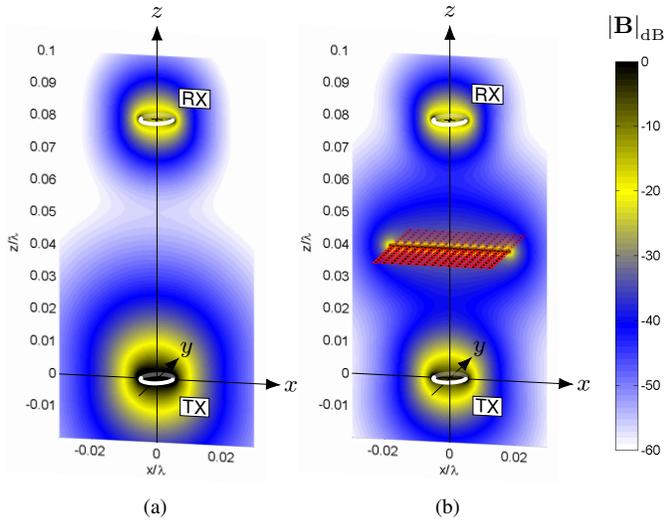

Fig. 10. The magnetic fields of the standard SISO system (a) and the SISO-P225 system (b), enhanced by an optimally loaded $15 \times 15$-element passive metasurface. The transmitter and receivers are located at a distance of $d = 0.08\lambda$, while the metasurface is placed in the middle, at $d/2$. The maximum achievable PTEs are about 1.2% for the reference (a) and 8.6% for the metasurface-enhanced case (b).

The graph in Fig. 9(b) confirms that the relaxation to obtain the maximum PTE with the passive metasurface (red solid) was tight in all cases under consideration, as the tightness error is in the range of $10^{-14}$ to $10^{-12}$.

Fig. 10 shows a comparison of the magnetic fields of the standard SISO WPT system and when it is optimally enhanced by the $15 \times 15$-element metasurface. Since the receiver load resistances differ very little and unit power is transferred in both cases, the fields around the receiver loops are very similar. However, the fields around the transmitter are quite different, since in the enhanced case the efficiency was improved from just above 1% to about 8.6%, requiring much lower input power and reducing the field amplitudes accordingly.

Lastly, in Fig. 11, the optimum load profiles of the metasurface elements are compared for three distances $d$. As can be seen, the closer the metasurface is to the transmitter and receiver, the larger the variance between different optimum capacitance values: the radial wave-like change in capacitance from the center of the metasurface to the outside at $0.01\lambda$ changes into a convex profile at $0.04\lambda$ and eventually into a concave profile at $0.08\lambda$ and further away. Edge effects due to the finite dimensions of the metasurface are responsible for the optimal values never being the same throughout the whole area.

## V. SUMMARY & CONCLUSION

The performance of WPT systems is limited by the size of the whole system. The requirement of high power transfer efficiency (PTE) generally puts a hard constraint on the distance between transmitter(s) and receiver(s); usually this distance can only be a small fraction of the wavelength at the operating frequency and is comparable to the transmitter loop size.

Many enhancements have been proposed to improve the performance of the WPT system. One of the most straightforward is to place passive elements in the proximity of the transmitter and/or receiver. Rigorously optimizing such systems is not a trivial task, since the required power constraints are nonconvex, thereby rendering the entire optimization problem unsolvable in general, at least for a large number of passive elements.

A fully equivalent work-around via tight semidefinite relaxation has been presented, which solves this optimization problem efficiently and reliably for an arbitrary number of passive elements. Starting from the impedance matrix of the whole WPT system and the load resistance of the receiver, the maximum achievable PTE as well as the required optimal currents and reactive loading elements are found. In an outer optimization loop, the load resistance can also be optimized; however it is generally observed that the optimum deviates very little from the one obtained from MISO optimization and that the maximum PTE is not sensitive to small deviations. A simple test of the solutions confirms tightness of the relaxation; i.e. that the solved problem remained equivalent to the original one and that the true global optimum was found.

With this powerful and versatile optimization method at hand, three examples of passively enhanced WPT systems have been investigated. It has been shown that optimally loaded passive elements can indeed lead to substantial improvement of the overall performance of such WPT systems. Relays can help mitigate problems of misalignment while also extending the region of high PTE. Multiple passive elements can be used cooperatively, to either form transmission-line type field guiding structures or more general structures such as metasurfaces and metamaterials.

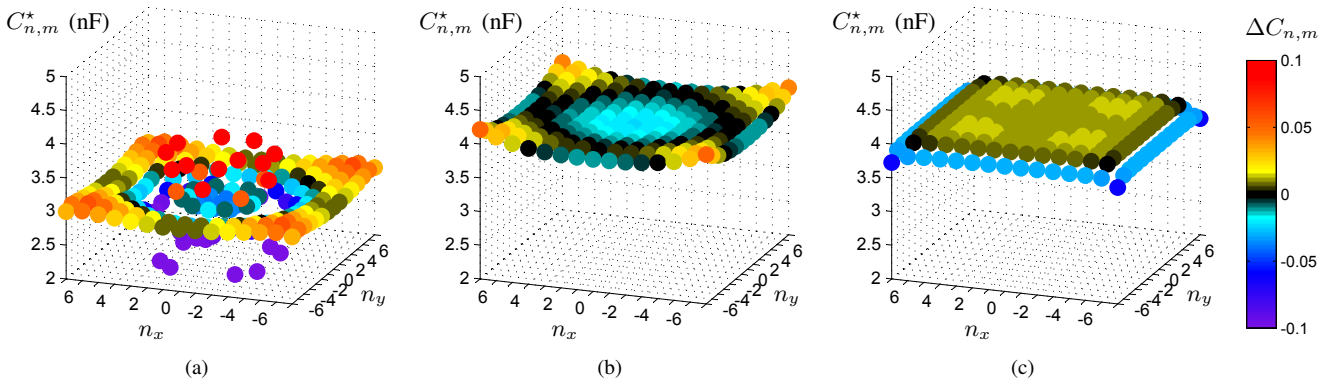

Fig. 11. Optimal capacitance profiles for the $15 \times 15$-element metasurface at distances (a) $d = 0.01\lambda = D_{\text{rx}}$, (b) $d = 0.04\lambda = 4D_{\text{rx}}$, (c) $d = 0.08\lambda = 8D_{\text{rx}}$.

# Appendix A
## Analytical Formulation of the Positive and Negative PIM Parts by Eigenvectors

Let the impedance matrix of the unloaded WPT system be

$$\mathbf{Z} = \mathbf{R} + j\omega\mathbf{L} + j\omega\mathbf{M} \,, \tag{26}$$

where $\mathbf{R}$ and $\mathbf{L}$ are real-valued diagonal matrices containing the losses and self-reactances (usually inductances $L_n > 0$, for loop-based magnetically coupled systems) of the transmitters and the receiver. $\mathbf{M} = \mathbf{M}^T$ is a symmetric (due to reciprocity of the passive system) hollow matrix containing the mutual impedances $j\omega M_{n,m}$. Generally, each $M_{n,m}$ is complex-valued, due to retardation effects when the electrical distance between the transmitter(s) and receiver are not very small.

Moreover, let

$$\{\lambda_{n,m}, \mathbf{v}_{n,m}\} = \mathrm{eig}(\mathbf{T}_n) \qquad n = 1, \dots, N \tag{27}$$

denote the $m$th eigenvalue $\lambda_{n,m}$ and corresponding eigenvector $\mathbf{v}_{n,m}$ of the $N \times N$ port impedance matrix (PIM) $\mathbf{T}_n$. In the order "negative, positive, and zero", the eigenvalues can be given as

$$\lambda_{n,m} = \begin{cases} -\lambda_n^- = \dfrac{1}{2}\left(R_n - \sqrt{S_n^2 + R_n^2}\right) < 0 & m = 1 \\[2mm] \lambda_n^+ = \dfrac{1}{2}\left(R_n + \sqrt{S_n^2 + R_n^2}\right) > 0 & m = 2 \\[2mm] 0 & m = 3, \dots, N \end{cases} \tag{28}$$

where the shorthand $S_n^2 = \omega^2 \sum_m |M_{n,m}|^2 > 0$ has been used. Note that, with these definitions, $\lambda_n^+, \lambda_n^- > 0$, for all $n$.

Finally, let the eigenvectors corresponding to the non-zero eigenvalues be denoted by $\mathbf{v}_{n,m=1,2} = \mathbf{v}_n^{\mp}$ (omitting the index $m$, as it is clear from (28) that $m = 1, 2$ correspond to the superscripts $-, +$, respectively) and a superscript zero point to eigenvectors corresponding to zero eigenvalues, i.e. $\mathbf{v}_{n,m>2} = \mathbf{v}_{n,m}^0$ (the index $m$ starts at 3, for these eigenvectors). Hence, the quadratic forms with respect to the eigenvalues and eigenvectors are given by

$$(\mathbf{v}_n^{\pm})^H \mathbf{T}_n \mathbf{v}_n^{\pm} = \pm\lambda_n^{\pm} \cdot (\mathbf{v}_n^{\pm})^H \mathbf{v}_n^{\pm} \tag{29a}$$

$$\mathbf{v}_{n,m}^H \mathbf{T}_n \mathbf{v}_{n,m} = 0 \qquad \forall \, m \neq 1, 2 \,. \tag{29b}$$

The PIMs can be separated into their positive and negative (semidefinite) parts

$$\mathbf{T}_n = \mathbf{T}_n^+ - \mathbf{T}_n^- \,, \tag{30}$$

where both $\mathbf{T}_n^+, \mathbf{T}_n^- \succeq 0$. Further, each part is simply obtained from its eigenvalues and eigenvectors

$$\mathbf{T}_n^{\pm} = \lambda_n^{\pm} \frac{\mathbf{v}_n^{\pm}(\mathbf{v}_n^{\pm})^H}{(\mathbf{v}_n^{\pm})^H \mathbf{v}_n^{\pm}} \tag{31}$$

with the denominator being the outer product and the numerator the inner product of the respective eigenvectors.

These eigenvectors $\mathbf{v}_n^{\pm}$, can be obtained analytically as well:

$$\mathbf{v}_n^{\pm} = \frac{1}{\tilde{M}_{n,N}^{\pm}} \begin{bmatrix} \tilde{M}_{n,1}^{\pm} \\ \vdots \\ \tilde{M}_{n,N}^{\pm} \end{bmatrix} \tag{32}$$

where $\tilde{M}_{n,m}^{\pm}$ are the entries of the matrix $\tilde{\mathbf{M}}$ which is identical to the mutual impedance matrix $\mathbf{M}$, as given in (26), with the exception of the diagonal:

$$\tilde{\mathbf{M}}^{\pm} = \mathbf{M} \pm 2\mathbf{Diag}(\lambda^{\pm}) \,. \tag{33}$$

Obtaining the positive and negative parts of the PIMs analytically, directly from the impedance matrix entries, adds both computational efficiency as well as numerical precision as compared to using a numerical eigenvalue decomposition.